\documentclass[10pt]{amsart} 
\RequirePackage[utf8]{inputenc}

\usepackage{amsmath,amsthm,amsfonts,amssymb,amscd,amsbsy,dsfont,hyperref}
\usepackage[all]{xy}

\usepackage{palatino}

\renewcommand{\geq}{\geqslant}
\renewcommand{\leq}{\leqslant}

\usepackage{geometry}
\allowdisplaybreaks

\hypersetup{
    pdftoolbar=true,
    pdfmenubar=true,
   pdffitwindow=false,
    pdfstartview={FitH},
    pdftitle={}, 
    pdfauthor={}, 
    pdfsubject={}, 
    pdfkeywords={},
    pdfnewwindow=true,
    colorlinks=true, 
    linkcolor=blue,
    citecolor=blue,
    urlcolor=black,}

\newtheorem*{acknowledgement}{Acknowledgement}
\newtheorem{corollary}{Corollary}

\newtheorem{theorem}{Theorem}

\numberwithin{equation}{section}

\title[Rigidity of 4-dimensional manifolds]{Rigidity of four-dimensional compact manifolds\\ with harmonic Weyl tensor}
\author[Ernani Ribeiro Jr]{Ernani Ribeiro Jr$^{1}$}

\address{Universidade Federal do Cear\'a - UFC, Departamento  de Matem\'atica, Campus do Pici, Av. Humberto Monte, Bloco 914,
60455-760-Fortaleza / CE, Brazil.}
 \email{ernani@mat.ufc.br$^{1}$}

\allowdisplaybreaks
\numberwithin{equation}{section}
\numberwithin{theorem}{section}

\thanks{$^{1}$Partially supported by grants from FUNCAP/Brazil and CNPq/Brazil.}
\date{September 21, 2014}

\begin{document}

\newcommand{\spacing}[1]{\renewcommand{\baselinestretch}{#1}\large\normalsize}
\spacing{1.2}

\begin{abstract}
The goal of this paper is to investigate the rigidity of 4-dimensional manifolds involving some pinching curvature conditions. To this end, we make use of the approach of biorthogonal curvature which is weaker than the sectional curvature. Here, we prove a rigi\-dity result for 4-dimensional compact manifolds under a suitable lower bound condition on the mi\-ni\-mum of the biorthogonal curvature. From this, we improve the pinching constants considered by some preceding works on a rigidity result for 4-dimensional manifolds.
 \end{abstract}

\maketitle

\section{Introduction}

It plays an important role in geometry to classify  4-dimensional compact
manifolds in the category of either topology, diffeomorphism, or isometry. This is because 4-dimensional  ma\-ni\-folds have special behavior. For instance, the bundle of $2$-forms can be invariantly decomposed as a direct sum; further relevant facts may be found in \cite{besse} and \cite{scorpan}. There has been a considerable amount of research on 4-dimensional manifolds involving some pinching curvature condition. We underline the next ones: nonnegative or positive sectional curvature (cf. \cite{Noronha} and \cite{Seaman}), nonnegative or positive Ricci curvature (cf. \cite{LNN} and \cite{Tani}), nonnegative or po\-si\-tive scalar curvature (cf. \cite{GX}, \cite{GX2} and \cite{Kotschick}), nonnegative or positive isotropic curvature (cf. \cite{brendle}, \cite{CTZ}, \cite{MW}, \cite{noronha2} and \cite{Seshadri}) and nonnegative or positive biorthogonal curvature (cf. \cite{renato}, \cite{CR} and \cite{Seaman}). Here, $M^4$ will denote a compact oriented 4-dimensional manifold and $g$ is a Riemannian metric on $M^4$ with scalar curvature $s_g,$ or simply $s,$ sectional curvature $K$ and volume form $dV_g.$

We now recall that  $(M^{4},\,g)$ is called {\it Einstein} if the Ricci curvature is given by $$Ric=\lambda g,$$ for some constant $\lambda.$ This means that $M^4$ has constant Ricci curvature. In a such case, a result due to Berger \cite{berger} combined with Synge's theorem  allows us to conclude that, if $M^4$ has positive sectional curvature, then it satisfies $$2\leq \chi(M)\leq  9,$$ where $\chi(M)$ stands for the Euler characteristic of $M^4.$ Furthermore, Gursky and LeBrun showed that $M^4$ satisfies the inequality $$\chi(M)\geq \frac{15}{4}|\tau(M)|,$$ where $\tau(M)$ denotes the signature of $M^4.$ This result improves the famous Hitchin-Thorpe inequality (cf. \cite{hitchin} and \cite{thorpe}). On the basis of these comments we may conclude that most 4-dimensional manifolds can not carry any Einstein structure with positive or nonnegative sectional curvature.

Tachibana \cite{tashibana} asserts that a compact Einstein manifold with positive curvature operator is isometric to a spherical space form. While Micallef and Wang \cite{MW} extended Tachibana's result for $n=4$ showing that a 4-dimensional compact Einstein manifold with nonnegative isotropic curvature is locally symmetric. Recently, Brendle \cite{brendle} proved that a compact Einstein manifold with positive isotropic curvature must be isometric to a spherical space form. 

In 2000, Yang \cite{Yang} showed a rigidity result for Eins\-tein structures with positive sectional curvature on  4-dimensional manifolds. More precisely, he proved the following result.

\begin{theorem}[Yang, \cite{Yang}]
\label{thmY}
Let $(M^4,\,g)$ be a 4-dimensional complete Einstein manifold with norma\-lized Ricci curvature $Ric=1.$ Suppose that
\begin{equation}
\label{1st}
K\geq \frac{\big(\sqrt{1249}-23\big)}{120}\approx 0.102843.
\end{equation} Then, $M^4$ is  isometric to either
\begin{enumerate}
\item $\Bbb{S}^4$ with its canonical
metric, or 
\item  $\Bbb{CP}^2$ with Fubini-Study metric.
\end{enumerate}
\end{theorem}

As it was pointed by Yang $0.102843$ is apparently not the best possible lower bound on the sectional curvature to get this conclusion. In fact, from a convergence argument it is possible to show that there is a constant $0<\beta< 0.102843$ such that the conclusion of Theorem \ref{thmY} even is true for $K\geq \beta.$ 

In 2004, Costa \cite{Costa} showed that Yang's result remains true under weaker condition 
\begin{equation}
\label{2nd}
K\geq\frac{\big(2-\sqrt{2}\big)}{6}\approx 0.09763.
\end{equation} 

It has been conjectured that:  

\begin{flushright}
\begin{minipage}[t]{4.37in}
 \emph{``Every 4-dimensional Einstein manifold with normalized Ricci curvature $Ric=1$ and positive sectional curvature must be isometric to either $\Bbb{S}^4$ or $\Bbb{CP}^2$ with their normalized standard metrics."}  \cite{Yang}.
 \end{minipage}
\end{flushright}

In order to explain our assumption in the main result to follow let us recall briefly the concept of biorthogonal curvature. For each plane $P\subset T_{p}M$ at a point $p\in M^4,$  we define  {\it the biorthogonal (sectional) curvature} of $P$ by the following average of the sectional curvatures
\begin{equation}\label{[1.2]}
 \displaystyle{K^\perp (P) = \frac{K(P) + K(P^\perp) }{2}},
\end{equation}
where $P^\perp$  is the orthogonal plane to  $P.$ For our purposes, for each point $p\in M^4,$ we take the minimum of biorthogonal curvature to obtain the following function
\begin{equation}
\label{[1.3]}
K_1^\perp(p) = \textmd{min} \{K^\perp(P); P  \textmd{ is a 2- plane in } T_{p}M \}.
\end{equation}

The sum of two sectional curvatures on two orthogonal planes  plays a very crucial role on 4-dimensional manifolds. This notion appeared previously in works due to Seaman \cite{Seaman} and Noronha \cite{Noronha}, see also \cite{LeBrun} (cf. Section 5). Interestingly enough, $\Bbb{S}^{1}\times \Bbb{S}^{3}$ with its canonical metric shows that the positivity of the biorthogonal curvature is an intermediate condition between positive sectional curvature and positive scalar curvature. It is not hard to prove that a 4-dimensional Riemanniana manifold $(M^4,\,g)$ is Einstein if and only if $K^\perp(P)=K(P)$ for any plane $P\subset T_{p}M$ at any point $p\in M^4$ (cf. Corollary 6.26 \cite{handbook} and \cite{ST}). For more details on this subject we address to \cite{renato}, \cite{CR}, \cite{Noronha}, \cite{noronha2} and \cite{Seaman}. 

Proceeding, we say that the Weyl tensor $W$ is harmonic if $\delta W = 0,$ where $\delta$ is the formal divergence defined for any
$(0,4)$-tensor $T$ by $$\delta T(X_1,X_2,X_3) =
-trace_{g}\{(Y,Z)\mapsto\nabla_{Y}T(Z,X_1,X_2,X_3)\},$$ where $g$ is
the metric of $M^4.$  We remark that metrics with harmonic curvature as well as conformally flat metrics with constant scalar curvature are real analytic in harmonic coordinates; see \cite{DeTurck}. One should be emphasized that every 4-dimensional Einstein manifold has harmonic Weyl tensor $W$ (cf. 16.65 in \cite{besse}, see also Lemma 6.14 in \cite{handbook}). Therefore, it is natural to ask which geometric implications has the assumption of the harmonicity of the tensor $W$ on a 4-dimensional manifold.

In 2013, inspired by Yang's work, Costa and the author proved the following result.

\begin{theorem}[Costa-Ribeiro Jr, \cite{CR}]
\label{thmT}
Let $(M^4,\, g)$  be a 4-dimensional compact oriented Riemannian manifold  with harmonic Weyl tensor and positive scalar curvature. We assume that $g$ is analytic  and  
\begin{equation}
\label{assump0}
K_1^{\perp}\geq\ \frac{s^{2}}{8(3s+5\lambda_{1})},
\end{equation} where $\lambda_{1}$ stands for  the first eigenvalue of the Laplacian with respect to $g.$ Then,  $M^4$ is either
\begin{enumerate}
\item  diffeomorphic to a connected sum  $\Bbb{S}^4 \sharp(\Bbb{R} \times \Bbb{S}^3)/G_1\sharp ...\sharp(\Bbb{R} \times \Bbb{S}^3)/G_n$, where each $G_i$ is a discrete subgroup of the isometry group of $\Bbb{R} \times \Bbb{S}^3.$ In this case, $g$ is locally conformally flat; or
\item isometric to a complex projective space  $\Bbb{CP}^2$ with the Fubini-Study metric.
\end{enumerate}
\end{theorem}

Inspired by the historical development on the study of the rigidity of 4-dimensional manifolds, in this paper, we use the notion of biorthogonal curvature to obtain a rigidity result for 4-dimensional compact manifold with harmonic Weyl tensor under a pinching condition weaker than (\ref{assump0}).

After these settings we may state our main result.

\begin{theorem}
\label{thmA}
Let $(M^4,\, g)$  be a 4-dimensional compact oriented Riemannian manifold  with harmonic Weyl tensor and positive scalar curvature. We assume that $g$ is analytic  and  
\begin{equation}
\label{assump}
K_1^{\perp}\geq\frac{s^2}{24(3\lambda_1+s)},
\end{equation} where $\lambda_{1}$ stands for  the first eigenvalue of the Laplacian with respect to $g.$ Then,  $M^4$ is either
\begin{enumerate}
\item  diffeomorphic to a connected sum  $\Bbb{S}^4 \sharp(\Bbb{R} \times \Bbb{S}^3)/G_1\sharp ...\sharp(\Bbb{R} \times \Bbb{S}^3)/G_n$, where each $G_i$ is a discrete subgroup of the isometry group of $\Bbb{R} \times \Bbb{S}^3.$ In this case, $g$ is locally conformally flat; or
\item isometric to a complex projective space  $\Bbb{CP}^2$ with the Fubini-Study metric.
\end{enumerate}
\end{theorem}

It is not difficult to check that $$\frac{s^2}{8(3s+5\lambda_1)}>\frac{s^2}{24(s+3\lambda_1)}.$$ For this, our condition (\ref{assump}) considered in Theorem \ref{thmA} is weaker that the former considered in Theorem \ref{thmT}. At same time, we already know that  Einstein metrics are analytic (cf. Theorem 5.26 in \cite{besse}). Also, as we have mentioned a 4-dimensional compact manifold $M^4$ is Einstein if and only if $K^\perp=K$. Furthermore, every Einstein metric has harmonic Weyl tensor. From these comments, we may conclude that Theorem \ref{thmA} generalizes Theorem \ref{thmY} and Theorem  \ref{thmT}.

One should be emphasized that the methods designed for the proof of Theorem \ref{thmA} are different from the proofs presented in \cite{CR} and \cite{Yang}. Our strategy is to use a more refined technique to assure a lower estimate to minimum of biorthogonal curvature.

In the sequel, it is not hard to prove that $Ric\geq\rho>0$ implies $\lambda_1\geq\frac{4\rho}{3}$ and $s\geq4\rho,$ see also \cite{gallot} and \cite{lich}. Then, we can combine Theorem \ref{thmA} with Tani's theorem \cite{Tani} and Bonnet-Myers theorem to obtain the following result.

\begin{corollary}\label{corRic}
Let $(M^4,g)$ be a 4-dimensional complete oriented Riemannian manifold with harmonic Weyl tensor and metric $g$ analytic. We
assume that $Ric\geq\rho>0$ and $K_1^\perp\geq\frac{s^2}{192\rho}.$ Then,  $M^4$ is  isometric to either
\begin{enumerate}
\item $\Bbb{S}^4$ with its canonical
metric, or 
\item $\Bbb{CP}^2$ with Fubini-Study metric.
\end{enumerate}
\end{corollary}

This improves Corollary 3 in \cite{CR}. Moreover, Corollary \ref{corRic} shows that the pinching constants used by Yang in (\ref{1st}) as well as by Costa in (\ref{2nd}) can be improved to $\approx 0.08333.$

\vspace{0.3cm}

The paper is organized as follows. In Section 2, we review some classical results on 4-dimensional manifolds that will be use here. Moreover, we briefly outline some useful informations on biorthogonal curvature. In Section 3, we prove the main result.

\section{Background}
Throughout this section we recall some informations and important results that will be useful in the proof of our main result.  In what follows $M^4$ will denote an oriented 4-dimensional manifold and $g$ is a Riemannian metric on $M^4.$ As it was previously pointed 4-manifolds are fairly special.  In fact, following the notations used in \cite{handbook} (see also \cite{VIA} p. 46), given any local orthogonal frame $\{e_{1}, e_{2}, e_{3}, e_{4}\}$ on open set of $M^4$ with dual basis $\{e^{1}, e^{2}, e^{3}, e^{4}\},$ there exists a unique bundle morphism $\ast$ called {\it Hodge star} (acting on bivectors) such that $$\ast(e^{1}\wedge e^{2})=e^{3}\wedge e^{4}.$$ With this setting we deduce that $\ast$ is an involution, i.e. $\ast^{2}=Id.$ In particular, this implies that the bundle of $2$-forms on a 4-dimensional oriented Riemannian manifold can be invariantly decomposed as a direct sum $$\Lambda^2=\Lambda^{2}_{+}\oplus\Lambda^{2}_{-}.$$ This allows us to conclude that the Weyl tensor $W$ is an endomorphism of $\Lambda^2=\Lambda^{+} \oplus \Lambda^{-} $ such that 
\begin{equation}
\label{decW}
W = W^+\oplus W^-.
\end{equation} 

Next, we recall that $dim_{\Bbb{R}}(\Lambda^2)=6$ and $dim_{\Bbb{R}}(\Lambda^{\pm})=3.$ Also, it is well-known that 
\begin{equation}
\label{6A}
\Lambda^{+}=span\Big\{\frac{e^{1}\wedge e^{2}+e^{3}\wedge e^{4}}{\sqrt{2}},\,\frac{e^{1}\wedge e^{3}+e^{4}\wedge e^{2}}{\sqrt{2}},\,\frac{e^{3}\wedge e^{2}+e^{4}\wedge e^{1}}{\sqrt{2}}\Big\}
\end{equation}
 and 
\begin{equation}
\label{6B}
 \Lambda^{-}=span\Big\{\frac{e^{1}\wedge e^{2}-e^{3}\wedge e^{4}}{\sqrt{2}},\,\frac{e^{1}\wedge e^{3}-e^{4}\wedge e^{2}}{\sqrt{2}},\,\frac{e^{3}\wedge e^{2}-e^{4}\wedge e^{1}}{\sqrt{2}}\Big\}.
\end{equation} Accordingly, the bundles $\Lambda^{+}$ and $\Lambda^{-}$ carry natural orientations such that the bases (\ref{6A}) and (\ref{6B}) are positive-oriented. For this, if $\mathcal{R}$ denotes the curvature of $M^4$ we get

\begin{equation}
\mathcal{R}=
\left(
  \begin{array}{c|c}
    \\
W^{+} +\frac{s}{12}Id & B \\ [0.4cm]\hline\\

    B^{*} & W^{-}+\frac{s}{12}Id  \\[0.4cm]
  \end{array}
\right),
\end{equation}
where $B:\Lambda^{-}\to \Lambda^{+}$ stands for the Ricci traceless operator of $M^4$ given by $B=Ric-\frac{s}{4}g.$ For more details see \cite{besse} and \cite{VIA}.

We now fix a point and diagonalize $W^\pm$ such that $w_i^\pm,$ $1\le i \le 3,$ are their respective eigenvalues. We stress that the eigenvalues of  $W^\pm$ satisfy 
\begin{equation}
\label{eigenvalues}
w_1^{\pm}\leq w_2^{\pm}\leq w_3^{\pm}\,\,\,\,\hbox{and}\,\,\,\,w_1^{\pm}+w_2^{\pm}+w_3^{\pm} = 0.
\end{equation} In particular, (\ref{eigenvalues}) allows us to infer
\begin{equation}
\label{estcos}
|W^\pm|^2\leq6(w_1^\pm)^2.
\end{equation} In fact, from (\ref{eigenvalues}) it is easy to see that $$(w_{2}^{\pm})^2+(w_{3}^{\pm})^2=(w_{1}^{\pm})^2-2w_{2}^{\pm}w_{3}^{\pm}.$$ Therefore, we have
\begin{eqnarray*}
|W^{+}|^2=2(w_{1}^{+})^2-2w_{2}^{+}w_{3}^{+}.
\end{eqnarray*} Now, taking into account that $ w_{1}^{+}w_{3}^{+}\leq w_{2}^{+}w_{3}^{+}$ and $(w_{1}^{+})^2=-w_{1}^{+}w_{3}^{+}-w_{1}^{+}w_{2}^{+}$ we deduce $|W^{+}|^2\leq6(w_{1}^{+})^2.$ In the same way we obtain $|W^{-}|^2\leq6(w_{1}^{-})^2.$

For the sake of completeness let us briefly outline the construction of the minimum of biorthogonal curvature; more details can be found in \cite{CR}, see also \cite{LeBrun} (cf. Section 5). To do so, we consider a point $p\in M^4$ and $X,Y\in T_{p}M$ orthonormal. Whence, the unitary 2-form $\alpha=X\wedge Y$ can be uniquely written as $\alpha=\alpha^{+}+\alpha^{-},$ where $\alpha^{\pm}\in \Lambda^{\pm}$ with $|\alpha^{+}|^{2}=\frac{1}{2}$ and $|\alpha^{-}|^{2}=\frac{1}{2}.$ From these settings, the sectional curvature $K(\alpha)$ can be written as
\begin{equation}
\label{k}
K(\alpha)=\frac{s}{12}+\langle\alpha^{+},W^{+}(\alpha^{+})\rangle+\langle\alpha^{-},W^{-}(\alpha^{-})\rangle+2\langle\alpha^{+},B\alpha^{-}\rangle.
\end{equation}
Moreover, we have
\begin{equation}
\label{kperp}
K(\alpha^{\perp})=\frac{s}{12}+\langle\alpha^{+},W^{+}(\alpha^{+})\rangle+\langle\alpha^{-},W^{-}(\alpha^{-})\rangle-2\langle\alpha^{+},B\alpha^{-}\rangle,
\end{equation} where $\alpha^{\perp}=\alpha^{+}-\alpha^{-}.$  Combining (\ref{k}) with (\ref{kperp}) we arrive at
\begin{equation}
\frac{K(\alpha)+K(\alpha^{\perp})}{2}=\frac{s}{12}+\langle\alpha^{+},W^{+}(\alpha^{+})\rangle+\langle\alpha^{-},W^{-}(\alpha^{-})\rangle.
\end{equation}
Hence, we may use (\ref{[1.3]}) to infer
\begin{eqnarray*}
\label{eq898}
K_{1}^{\perp}=\frac{s}{12}+min\left\{\langle\alpha^{+},W^{+}(\alpha^{+})\rangle; \,|\alpha^{+}|^{2}=\frac{1}{2}\right\}+min\left\{\langle\alpha^{-},W^{-}(\alpha^{-})\rangle; \,|\alpha^{-}|^{2}=\frac{1}{2}\right\}.
\end{eqnarray*}
Furthermore, as it was explained in \cite{CR} and \cite{Seaman}, Equation (\ref{[1.3]}) provides the following useful identity
\begin{equation}
\label{[1.6]}
K_1^\perp = \frac{w_1^+ + w_1^-}{2}+\frac{s}{12}.
\end{equation}

Proceeding, given a section $T\in\Gamma(E),$ where
$E\to M$ is a vector bundle over $M,$ we already know that  $$|\nabla|T||\leq|\nabla T|$$ away from the zero locus of $T.$  This inequality is known as  {\it Kato's inequality}. In \cite{gursky2}, Gursky and LeBrun proved a {\it refined Kato's inequality}. More precisely, if $W^+$ is harmonic, then away from the zero locus of $W^+$ we have
\begin{equation}
\label{refkato} 
|\nabla|W^{+}||\leq \sqrt{\frac{3}{5}}|\nabla W^{+}|,
\end{equation} for more details see Lemma 2.1 in \cite{gursky}.

We also recall that if $\delta W^\pm=0,$ then the Weitzenb\"{o}ck formula (cf. 16.73 in \cite{besse}, see also \cite{handbook}) is given by
\begin{equation}
\label{weitzenbock}
\Delta |W^{\pm}|^2=2|\nabla W^{\pm}|^2+s|W^{\pm}|^2- 36 \det W^{\pm}.
\end{equation} We highlight that our definition of Laplacian differs from \cite{besse} by sign, i.e. $\Delta f={\rm tr}Hess\, f.$ Finally, by using Lagrange-multipliers it is not difficult to check the following inequality
\begin{equation}
\label{detest}
\det W^{+}\leq\frac{\sqrt{6}}{18}|W^{+}|^{3}.
\end{equation}

\section{Proof of the Main Result}

\subsection{Proof of Theorem \ref{thmA}}

\begin{proof}
To begin with, we assume that $(M^4,\, g)$ is a compact oriented Riemannian manifold with positive scalar curvature and harmonic Weyl tensor. Since $g$ is analytic we deduce that $|W^{\pm}|^2$ are analytic. Moreover, supposing $W^{\pm}\not\equiv 0$ it is easy to see that the set $$\Sigma=\Big\{p\in M;|W^+|(p)=0 \mbox{ or }|W^-|(p)=0\Big\}$$ is finite. 

We now suppose by contradiction that $(M^4,g)$ is not half conformally flat. For this, for some $\alpha>0$ (to be chosen later) and any $\varepsilon>0,$ there exists $t=t(\alpha,\varepsilon)>0$ such that 
$$\int_M\Big(\left(|W^+|^2+\varepsilon\right)^{\frac{\alpha}{2}}-t\left(|W^-|^2+\varepsilon\right)^{\frac{\alpha}{2}}\Big)dV_g=0.$$

On the other hand, we notice that 

\begin{eqnarray*}
\label{eq1a}
\Delta \Big(\big(|W^{+}|^2 +\varepsilon\big)^{\alpha}&+&t^{2}\big(|W^{-}|^2 +\varepsilon\big)^{\alpha}\Big)\\&&=\alpha\big(|W^{+}|^2 +\varepsilon\big)^{\alpha-2}\Big(\big(|W^{+}|^{2}+\varepsilon\big)\Delta |W^{+}|^{2}+(\alpha-1)|\nabla|W^{+}|^2|^2\Big)\nonumber\\&&+t^{2}\alpha\big(|W^{-}|^2 +\varepsilon\big)^{\alpha-2}\Big(\big(|W^{-}|^{2}+\varepsilon\big)\Delta |W^{-}|^{2}+(\alpha-1)|\nabla|W^{-}|^2|^2\Big).
\end{eqnarray*} We recall that in dimension 4 we have $$|\delta W|^2=|\delta W^{+}|^2+|\delta W^{-}|^2.$$ Therefore, since $\delta W^\pm=0$ we invoke the Weitzenb\"{o}ck formula (\ref{weitzenbock}) to arrive at

\begin{eqnarray*}
\label{eq1a}
\Delta \Big(\big(|W^{+}|^2 &+&\varepsilon\big)^{\alpha}+t^{2}\big(|W^{-}|^2 +\varepsilon\big)^{\alpha}\Big)\\&=&\alpha\big(|W^{+}|^2 +\varepsilon\big)^{\alpha-2}\Big(\big(|W^{+}|^{2}+\varepsilon\big)\big(2|\nabla W^{+}|^2 +s|W^{+}|^2 -36\det W^{+}\big)\nonumber\\&&+(\alpha-1)|\nabla|W^{+}|^2|^2\Big)+t^{2}\alpha\big(|W^{-}|^2 +\varepsilon\big)^{\alpha-2}\Big(\big(|W^{-}|^{2}+\varepsilon\big)\big(2|\nabla W^{-}|^2+s|W^{-}|^2\nonumber\\&&-36\det W^{-}\big)+(\alpha-1)|\nabla|W^{-}|^2|^2\Big).
\end{eqnarray*} Upon integrating of the above expression over $M^4$ we obtain

\begin{eqnarray}\label{pteq1}
0&=&\alpha\int_M\left(|W^+|^2+\varepsilon\right)^{\alpha-1}(s|W^+|^2-36\det W^+)dV_g\nonumber\\
 &&+\alpha t^2\int_M\left(|W^-|^2+\varepsilon\right)^{\alpha-1}(s|W^-|^2-36\det W^-)dV_g\nonumber\\
 &&+\alpha\int_M\left(|W^+|^2+\varepsilon\right)^{\alpha-2}\Big(2\left(|W^+|^2+\varepsilon\right)|\nabla W^+|^2+(\alpha-1)|\nabla|W^+|^2|^2\Big)dV_g\nonumber\\
 &&+\alpha t^2\int_M\left(|W^-|^2+\varepsilon\right)^{\alpha-2}\Big(2\left(|W^-|^2+\varepsilon\right)|\nabla W^-|^2+(\alpha-1)|\nabla|W^-|^2|^2\Big)dV_g.
\end{eqnarray}

Next, we apply the refined Kato's inequality (\ref{refkato}) as well as (\ref{detest}) in (\ref{pteq1}) to infer

\begin{eqnarray*}
0&\geq&\alpha\int_M\left(|W^+|^2+\varepsilon\right)^{\alpha-1}(s|W^+|^2-2\sqrt{6}|W^+|^3)dV_g\\
 &&+\alpha t^2\int_M\left(|W^-|^2+\varepsilon\right)^{\alpha-1}(s|W^-|^2-2\sqrt{6}|W^-|^3)dV_g\\
 &&+\alpha\int_M\left(|W^+|^2+\varepsilon\right)^{\alpha-2}\Big(\frac{10}{3}|W^+|^2|\nabla|W^+||^2+\frac{10}{3}\varepsilon|\nabla|W^+||^2+(\alpha-1)|\nabla|W^+|^2|^2\Big)dV_g\\ &&+\alpha t^2\int_M\left(|W^-|^2+\varepsilon\right)^{\alpha-2}\Big(\frac{10}{3}|W^-|^2|\nabla|W^-||^2+\frac{10}{3}\varepsilon|\nabla|W^-||^2+(\alpha-1)|\nabla|W^-|^2|^2\Big)dV_g.
\end{eqnarray*}

Easily one verifies that 

\begin{equation*}
\frac{1}{4}|\nabla|W^{\pm}|^2|^2=|W^{\pm}|^2|\nabla|W^{\pm}||^2.
\end{equation*} 

From this, we deduce

\begin{eqnarray}
\label{p12}
0&\geq&\alpha\int_M|W^+|^2\left(|W^+|^2+\varepsilon\right)^{\alpha-1}(s-2\sqrt{6}|W^+|)dV_g\nonumber\\
 &&+\alpha t^2\int_M|W^-|^2\left(|W^-|^2+\varepsilon\right)^{\alpha-1}(s-2\sqrt{6}|W^-|)dV_g\nonumber\\
 &&+\alpha\int_M\left(|W^+|^2+\varepsilon\right)^{\alpha-2}\Big(\frac{5}{6}|\nabla|W^+|^2|^2+(\alpha-1)|\nabla|W^+|^2|^2\Big)dV_g\nonumber\\
 &&+\alpha t^2\int_M\left(|W^-|^2+\varepsilon\right)^{\alpha-2}\Big(\frac{5}{6}|\nabla|W^-|^2|^2+(\alpha-1)|\nabla|W^-|^2|^2\Big)dV_g\nonumber\\
 &=&\alpha\int_M|W^+|^2\left(|W^+|^2+\varepsilon\right)^{\alpha-1}(s-2\sqrt{6}|W^+|)dV_g\nonumber\\
 &&+\alpha t^2\int_M|W^-|^2\left(|W^-|^2+\varepsilon\right)^{\alpha-1}(s-2\sqrt{6}|W^-|)dV_g\nonumber\\
 &&+\alpha\left(\alpha-\frac{1}{6}\right)\int_M\Big(\left(|W^+|^2+\varepsilon\right)^{\alpha-2}|\nabla|W^+|^2|^2+t^2\left(|W^-|^2+\varepsilon\right)^{\alpha-2}|\nabla|W^-|^2|^2\Big) dV_g.
\end{eqnarray} Moreover, a straightforward  computation gives $$\alpha\left(|W^\pm|^2+\varepsilon\right)^{\alpha-2}|\nabla|W^\pm|^2|^2=\frac{4}{\alpha}|\nabla\left(|W^\pm|^2+\varepsilon\right)^{\frac{\alpha}{2}}|^2.$$ This jointly with (\ref{p12}) yields 
\begin{eqnarray}
\label{ptin1}
0&\geq&\alpha\int_M|W^+|^2\left(|W^+|^2+\varepsilon\right)^{\alpha-1}(s-2\sqrt{6}|W^+|)dV_g\nonumber\\
 &&+\alpha t^2\int_M|W^-|^2\left(|W^-|^2+\varepsilon\right)^{\alpha-1}(s-2\sqrt{6}|W^-|)dV_g\\
 &&+\left(4-\frac{2}{3\alpha}\right)\int_M\Big(|\nabla\left(|W^+|^2+\varepsilon\right)^{\frac{\alpha}{2}}|^2+t^2|\nabla\left(|W^-|^2+\varepsilon\right)^{\frac{\alpha}{2}}|^2\Big)dV_g.\nonumber
\end{eqnarray}

On the other hand, we have
\begin{eqnarray}
\label{ptin2}
|\nabla\left(|W^+|^2+\varepsilon\right)^{\frac{\alpha}{2}}|^2+t^2|\nabla\left(|W^-|^2+\varepsilon\right)^{\frac{\alpha}{2}}|^2&=&\frac{1}{2}\left\{|\nabla\big(\left(|W^+|^2+\varepsilon\right)^{\frac{\alpha}{2}}-t\left(|W^-|^2+\varepsilon\right)^{\frac{\alpha}{2}}\big)|^2\right.\nonumber\\
 &&+\left.|\nabla\big(\left(|W^+|^2+\varepsilon\right)^{\frac{\alpha}{2}}+t\left(|W^-|^2+\varepsilon\right)^{\frac{\alpha}{2}}\big)|^2\right\}\nonumber\\
 &\geq&\frac{1}{2}|\nabla\big(\left(|W^+|^2+\varepsilon\right)^{\frac{\alpha}{2}}-t\left(|W^-|^2+\varepsilon\right)^{\frac{\alpha}{2}}\big)|^2.
\end{eqnarray} At same time, from Poincar\'e inequality we get
\begin{equation}
\label{ptin3}
\displaystyle{\int_M|\nabla\big(\left(|W^+|^2+\varepsilon\right)^{\frac{\alpha}{2}}-t\left(|W^-|^2+\varepsilon\right)^{\frac{\alpha}{2}}\big)|^2 dV_{g}\geq\lambda_1\int_M\Big(\left(|W^+|^2+\varepsilon\right)^{\frac{\alpha}{2}}-t\left(|W^-|^2+\varepsilon\right)^{\frac{\alpha}{2}}\Big)^2dV_g}.
\end{equation} So, we combine (\ref{ptin2}) with (\ref{ptin3}) to infer

\begin{equation}
\label{pq1}
\displaystyle{\int_M |\nabla\left(|W^+|^2+\varepsilon\right)^{\frac{\alpha}{2}}|^2+t^2|\nabla\left(|W^-|^2+\varepsilon\right)^{\frac{\alpha}{2}}|^2 dV_g\geq\frac{\lambda_1}{2}\int_M\Big(\left(|W^+|^2+\varepsilon\right)^{\frac{\alpha}{2}}-t\left(|W^-|^2+\varepsilon\right)^{\frac{\alpha}{2}}\Big)^2dV_g}.
\end{equation}

Hereafter, we compare (\ref{pq1}) with (\ref{ptin1}) and pick $\alpha$ such that $\big(4-\frac{2}{3\alpha}\big)\geq0$ to arrive at
\begin{eqnarray*}
0&\geq&\int_M|W^+|^2\left(|W^+|^2+\varepsilon\right)^{\alpha-1}(s-2\sqrt{6}|W^+|)dV_{g}\nonumber\\&&+\int_{M}t^2|W^-|^2\left(|W^-|^2+\varepsilon\right)^{\alpha-1}(s-2\sqrt{6}|W^-|)dV_g\nonumber\\&&+\frac{1}{\alpha}\left(2-\frac{1}{3\alpha}\right)\lambda_1\Big[\int_M\Big(\left(|W^+|^2+\varepsilon\right)^\alpha-2t\left(|W^+|^2+\varepsilon\right)^{\frac{\alpha}{2}}\left(|W^-|^2+\varepsilon\right)^{\frac{\alpha}{2}}\Big)dV_g\nonumber\\&&+\int_{M}t^2\left(|W^-|^2+\varepsilon\right)^\alpha dV_{g}\Big].
\end{eqnarray*} Moreover, choosing $\alpha_0=\frac{1}{3},$ which maximizes $\frac{1}{\alpha}\left(2-\frac{1}{3\alpha}\right),$ we get
\begin{eqnarray*}
0&\geq&\int_M\Big(|W^+|^2\left(|W^+|^2+\varepsilon\right)^{\alpha_0-1}(s-2\sqrt{6}|W^+|)\Big)dV_g\nonumber\\&&+\int_{M}t^2|W^-|^2\left(|W^-|^2+\varepsilon\right)^{\alpha_0-1}(s-2\sqrt{6}|W^-|)dV_g+3\lambda_1\Big[\int_{M}\left(|W^+|^2+\varepsilon\right)^{\alpha_0}dV_{g}\nonumber\\&&-2\int_{M}t\left(|W^+|^2+\varepsilon\right)^{\frac{\alpha_0}{2}}\left(|W^-|^2+\varepsilon\right)^{\frac{\alpha_0}{2}}dV_{g}+\int_{M}t^2\left(|W^-|^2+\varepsilon\right)^{\alpha_0}dV_g\Big].
\end{eqnarray*} Therefore, when $\varepsilon$ goes to $0$ we obtain

\begin{eqnarray}
\label{ptin4}
0&\geq&\int_M\Big(|W^-|^{2\alpha_0}(s+3\lambda_1-2\sqrt{6}|W^-|)t^2-6\lambda_1|W^+|^{\alpha_0}|W^-|^{\alpha_0}t\nonumber\\
 &&+|W^+|^{2\alpha_0}(s+3\lambda_1-2\sqrt{6}|W^+|)\Big)dV_g.
\end{eqnarray}

We now remark that the integrand of (\ref{ptin4}) is a quadratic function of $t.$ Whence, for sim\-pli\-city, we set $$P(t)=|W^-|^{2\alpha_0}(a-2\sqrt{6}|W^-|)t^2-6\lambda_1|W^+|^{\alpha_0}|W^-|^{\alpha_0}t+|W^+|^{2\alpha_0}(a-2\sqrt{6}|W^+|),$$ where $a=s+3\lambda_1.$ 

We also notice that the discriminant $\Delta$ of $P(t)$ is given by
\begin{equation}
\label{ptdisc}\Delta=36\lambda_1^2|W^+|^{2\alpha_0}|W^-|^{2\alpha_0}-4|W^-|^{2\alpha_0}|W^+|^{2\alpha_0}(a-2\sqrt{6}|W^-|)(a-2\sqrt{6}|W^+|).
\end{equation}

In the sequel we claim that $\Delta$ is less than or equal to zero. In fact, from (\ref{estcos}) we already know that $$|W^\pm|^2\leq6(w_1^\pm)^2.$$ This together with (\ref{[1.6]}) provides
\begin{equation}
\label{ptin5}|W^+|+|W^-|\leq\sqrt{6}\Big(\frac{s}{6}-2K_1^\perp\Big).
\end{equation} Moreover, a straightforward computation using our assumption on $K_1^\perp$ yields
\begin{equation}
\label{ptin6}\left(\frac{s}{6}-2K_1^\perp\right)\leq\frac{a^2-9\lambda_1^2}{12a}.
\end{equation}
Next, we combine (\ref{ptin5}) with (\ref{ptin6}) and (\ref{ptdisc}) to deduce
\begin{eqnarray*}
\Delta&=&|W^-|^{2\alpha_0}|W^+|^{2\alpha_0}\Big(36\lambda_1^2-4a^2+8\sqrt{6}a(|W^+|+|W^-|)-96|W^+||W^-|\Big)\\
 &\leq&|W^-|^{2\alpha_0}|W^+|^{2\alpha_0}\Big(36\lambda_1^2-4a^2+4(a^2-9\lambda_1^2)-96|W^+||W^-|\Big)\\
 &=&-96|W^+|^{2\alpha_{0}+1}|W^-|^{2\alpha_{0}+1}\\
 &\leq&0,
\end{eqnarray*} which settles our claim. Hereafter, we use once more  (\ref{ptin4}) to conclude $|W^{+}||W^{-}|=0$ in $M^4$. From this, since  $\Sigma$ is finite we arrive at a contradiction. Therefore, we have shown that $(M^4,\,g)$ is half conformally flat.

From now on it suffices to follow the arguments applied in the final steps of the the proof of Theorem 6 in \cite{CR} (see also \cite{Yang}). More precisely, we define the following set $$A=\Big\{p\in M^{4};\,Ric(p)\neq\frac{s(p)}{4}g\Big\},$$ where $(Ric-\frac{s}{4}g)$ stands for the traceless Ricci tensor of $(M^{4},\,g).$ Then, if $A$ is empty we use Hitchin's theorem \cite{hitchin} to deduce that $M^4$ is either is isometric to $\Bbb{S}^4$ with its canonical metric or isometric to $\Bbb{CP}^{2}$ with Fubini-Study metric. Otherwise, if $A$ is not empty we deduce that $M^4$ is locally conformally flat. In other words, one the following assertions holds:
\begin{enumerate}
\item $M^4$ is isometric to $\Bbb{S}^4$ with its canonical metric;
\item $M^4$ isometric to $\Bbb{CP}^{2}$ with Fubini-Study metric; 
\item or $M^4$ has positive isotropic curvature.
\end{enumerate}
In this last case we can invoke Chen-Tang-Zhu theorem  \cite{CTZ} to conclude that $M^4$ is diffeomorphic to a connected sum  $\Bbb{S}^4 \sharp(\Bbb{R} \times \Bbb{S}^3)/G_1\sharp ...\sharp(\Bbb{R} \times \Bbb{S}^3)/G_n$, where each $G_i$ is a discrete subgroup of the isometry group of $\Bbb{R} \times \Bbb{S}^3$ (see also \cite{MW}). 

So, we finish the proof of Theorem \ref{thmA}.
\end{proof}

\begin{acknowledgement}
The author would like to thank the Depart. of Mathematics - Lehigh University, where this work was carried out. He is grateful to H.-D. Cao for the warm hospitality and his constant support. Moreover, he would like to extend his special thanks to E. Costa, P. Wu and R. Di\'ogenes for many helpful conversations  that benefited the presentation of this paper.
\end{acknowledgement}

\end{document}